\documentclass{amsart}
\usepackage{%
  mathtools,  
  amssymb,    
  amsthm,     
  amsfonts,   
  thmtools,   
  graphicx,
  float,
  color,
  xcolor,
  tikz,
  tikz-cd,
  mathrsfs,   
  etoolbox,   
  quiver
}

\definecolor{darkred}{rgb}{0.75,0,0}
\def\customcitecolor{darkred}
\def\customlinkcolor{darkred}

\usepackage[%
    colorlinks,
    citecolor=\customcitecolor,%
    linkcolor=\customlinkcolor,%
    urlcolor=\customlinkcolor%
]{hyperref}

\usepackage[capitalise,nameinlink,noabbrev]{cleveref}

\usepackage[margin=1in]{geometry}

\usepackage[final,nopatch=footnote]{microtype}

\theoremstyle{definition}
\newtheorem{theorem}{Theorem}[section]
\numberwithin{theorem}{section} 
\numberwithin{equation}{section} 

\usepackage{mfirstuc}
\newcommand{\capitalizename}[1]{\makefirstuc{#1}}

\newcommand{\defthm}[1]{%
  \newtheorem{#1}[theorem]{\capitalizename{#1}}%
}

\newcommand{\defthms}[1]{%
  \forcsvlist{\defthm}{#1}%
}

\defthms{%
  answer,assumption,claim,computation,conjecture,construction,corollary,
  counterexample,definition,digression,discussion,example,
  examples,exercise,fact,goal,idea,intuition,lemma,
  motivation,notation,note,proposition,question,remark,setup,
  slogan,strategy,terminology,upshot,warning%
}

\Crefname{computation}{Computation}{Computations}
\Crefname{question}{Question}{Questions}
\Crefname{conjecture}{Conjecture}{Conjectures}

\newcommand{\deftextcommand}[1]{%
  \expandafter\providecommand\csname #1\endcsname{\mathrm{#1}}%
}
\newcommand{\deftextcommands}[1]{%
    \forcsvlist{\deftextcommand}{#1}%
}

\deftextcommands{ab,alg,an,ann,Aut,BG,BGL,Bl,BO,BP,BSL,BSO,BSp,BSU,BU,can,cd,cdh,cl,coBar,codim,codom,coeq,coev,cof,cofib,coker,colim,coim,cone,conj,const,coTor,cyc,diag,Desc,dg,Disc,disc,dual,eff,EKL,End,eq,ess,et,Et,EU,ev,Ex,ex,Exc,Ext,fib,Fix,Fl,FP,fppf,fpqc,Frac,Frob,Fun,Gal,gen,GL,gp,Gr,gr,GW,Her,Ho,hocofib,hocolim,hofib,holim,Hom,id,Idem,im,incl,Ind,ind,inj,Inn,Inv,inv,iso,Jac,KGL,kgl,KH,KO,ko,KQ,kq,KR,KSp,KU,ku,Lan,Map,map,MGL,MO,Mor,mor,MSL,MSO,MSp,MSU,MU,mult,MUP,Nm,ob,obj,op,Orb,ord,Out,perf,Perm,PGL,pr,pre,Proj,proj,prom,PSL,quot,Ran,rank,Res,res,RO,sep,sgn,SH,sig,Sing,SL,SO,soc,Sp,Span,Spin,spn,Sq,st,Stab,SU,supp,Supp,Syl,syl,Sym,syn,SYT,TC,td,Th,THH,Tor,Tot,TP,TR,Tr,tr,triv,univ,var,veff,vol,Wel,Wr}

\deftextcommands{Ab,Aff,Alg,Ani,Bimod,CAlg,Cat,CDGA,CG,CGWH,Ch,CMon,coAlg,Coh,CommRing,ConjSub,coMod,Cor,Corr,CoSh,CRing,CW,Field,Fin,FinSet,Gpd,Grp,Grpd,Grph,Kan,Kar,LMod,Mack,Mfld,Mod,NAlg,Ouv,Perf,Poset,Pr,Pre,PSh,PShv,qCat,QCoh,Rep,Ring,RMod,sAb,Set,SH,Sh,Shv,Sm,Sp,Spc,Spectra,sPre,sSet,sShv,Stack,Sub,Top,Tors,Var,Vect}

\newcommand{\defblackboardletter}[1]{%
  \expandafter\providecommand\csname #1\endcsname{\mathbb{#1}}
}
\newcommand{\defblackboardletters}[1]{%
  \forcsvlist{\defblackboardletter}{#1}%
}

\defblackboardletters{A,C,F,P,Q,R,Z}

\providecommand{\xto}[1]{\xrightarrow{#1}}

\RequirePackage{bbm}

\RequirePackage{pict2e,picture}
\makeatletter
\DeclareRobustCommand{\DDelta}{{\mathpalette\bb@Delta\relax}}
\newcommand{\bb@Delta}[2]{%
  \begingroup
  \sbox\z@{$\m@th#1\Delta$}%
  \dimendef\Dht=6 \dimendef\Dwd=8
  \setlength{\Dwd}{\wd\z@}%
  \setlength{\Dht}{\ht\z@}%
  \begin{picture}(\Dwd,\Dht)
  \put(0,0){$\m@th#1\Delta$}
  \put(.42\Dwd,.7\Dht){\line(10,-26){.25\Dht}}
  \end{picture}%
  \endgroup
}

\usepackage{listings}
\usepackage{xcolor}
\definecolor{darkred}{rgb}{0.75,0,0}
\lstdefinelanguage{Macaulay2}
{
basicstyle={\ttfamily},
keywordstyle={\color{darkred!80!black}},
commentstyle={\color{gray}},
stringstyle={\color{red!40!black}},
rulecolor=\color{darkred},
basewidth={1.2ex}, 
sensitive=false,
morecomment=[l]{--},
morecomment=[s]{-*}{*-},
morestring=[b]",
escapechar={`},
escapebegin={\rmfamily},
morekeywords={about,abs,AbstractToricVarieties,accumulate,Acknowledgement,acos,acosh,acot,addCancelTask,addDependencyTask,addEndFunction,addHook,AdditionalPaths,addStartFunction,addStartTask,Adjacent,adjoint,AdjointIdeal,AffineVariety,AfterEval,AfterNoPrint,AfterPrint,agm,AInfinity,alarm,AlgebraicSplines,Algorithm,Alignment,all,AllCodimensions,allowableThreads,ambient,analyticSpread,Analyzer,AnalyzeSheafOnP1,ancestor,ancestors,ANCHOR,and,andP,AngleBarList,ann,annihilator,antipode,any,append,applicationDirectory,applicationDirectorySuffix,apply,applyKeys,applyPairs,applyTable,applyValues,apropos,argument,Array,arXiv,Ascending,ascii,asin,asinh,ass,assert,associatedGradedRing,associatedPrimes,AssociativeAlgebras,AssociativeExpression,atan,atan2,atEndOfFile,Authors,autoload,AuxiliaryFiles,backtrace,Bag,Bareiss,baseFilename,BaseFunction,baseName,baseRing,baseRings,BaseRow,BasicList,basis,BasisElementLimit,Bayer,BeforePrint,beginDocumentation,BeginningMacaulay2,Benchmark,benchmark,Bertini,BesselJ,BesselY,betti,BettiCharacters,BettiTally,between,BGG,BIBasis,Binary,BinaryOperation,Binomial,binomial,BinomialEdgeIdeals,Binomials,BKZ,BlockMatrix,BLOCKQUOTE,BODY,Body,BoijSoederberg,BOLD,Book3264Examples,Boolean,BooleanGB,borel,Boxes,BR,break,Browse,Bruns,cache,CacheExampleOutput,CacheFunction,CacheTable,cacheValue,CallLimit,cancelTask,capture,catch,Caveat,CC,CDATA,ceiling,Center,centerString,Certification,ChainComplex,chainComplex,ChainComplexExtras,ChainComplexMap,ChainComplexOperations,ChangeMatrix,char,CharacteristicClasses,characters,charAnalyzer,check,CheckDocumentation,chi,Chordal,class,Classic,clean,clearAll,clearEcho,clearOutput,close,closeIn,closeOut,ClosestFit,CODE,code,codim,CodimensionLimit,coefficient,CoefficientRing,coefficientRing,coefficients,Cofactor,CohenEngine,CohenTopLevel,CoherentSheaf,CohomCalg,cohomology,coimage,CoincidentRootLoci,coker,cokernel,collectGarbage,columnAdd,columnate,columnMult,columnPermute,columnRankProfile,columnSwap,combine,Command,commandInterpreter,commandLine,COMMENT,commonest,commonRing,comodule,CompactMatrix,compactMatrixForm,CompiledFunction,CompiledFunctionBody,CompiledFunctionClosure,Complement,complement,complete,CompleteIntersection,CompleteIntersectionResolutions,Complexes,ComplexField,components,compose,compositions,compress,concatenate,conductor,ConductorElement,cone,Configuration,ConformalBlocks,conjugate,connectionCount,Consequences,Constant,Constants,constParser,content,continue,contract,Contributors,ConvexInterface,conwayPolynomial,ConwayPolynomials,copy,copyDirectory,copyFile,copyright,Core,CorrespondenceScrolls,cos,cosh,cot,CotangentSchubert,cotangentSheaf,coth,cover,coverMap,cpuTime,createTask,Cremona,csc,csch,current,currentColumnNumber,currentDirectory,currentFileDirectory,currentFileName,currentLayout,currentLineNumber,currentPackage,currentString,currentTime,Cyclotomic,Database,Date,DD,dd,deadParser,debug,debugError,DebuggingMode,debuggingMode,debugLevel,DecomposableSparseSystems,Decompose,decompose,deepSplice,Default,default,defaultPrecision,Degree,degree,degreeLength,DegreeLift,DegreeLimit,DegreeMap,DegreeOrder,DegreeRank,Degrees,degrees,degreesMonoid,degreesRing,delete,demark,denominator,Dense,Density,Depth,depth,Descending,Descent,Describe,describe,Description,det,determinant,DeterminantalRepresentations,DGAlgebras,diagonalMatrix,diameter,Dictionary,dictionary,dictionaryPath,diff,DiffAlg,difference,dim,directSum,disassemble,discriminant,dismiss,Dispatch,distinguished,DIV,Divide,divideByVariable,DivideConquer,DividedPowers,Divisor,DL,Dmodules,do,doc,docExample,docTemplate,document,DocumentTag,Down,drop,DT,dual,eagonNorthcott,EagonResolution,echoOff,echoOn,EdgeIdeals,edit,EigenSolver,eigenvalues,eigenvectors,eint,EisenbudHunekeVasconcelos,elapsedTime,elapsedTiming,elements,Eliminate,eliminate,Elimination,EliminationMatrices,EllipticCurves,EllipticIntegrals,else,EM,Email,End,end,endl,endPackage,Engine,engineDebugLevel,EngineRing,EngineTests,entries,EnumerationCurves,environment,Equation,EquivariantGB,erase,erf,erfc,error,errorDepth,euler,EulerConstant,eulers,even,EXAMPLE,ExampleFiles,ExampleItem,examples,ExampleSystems,Exclude,exec,exit,exp,expectedReesIdeal,expm1,exponents,export,exportFrom,exportMutable,Expression,expression,Ext,extend,ExteriorIdeals,ExteriorModules,exteriorPower,Factor,factor,false,Fano,FastMinors,FastNonminimal,FGLM,File,fileDictionaries,fileExecutable,fileExists,fileExitHooks,fileLength,fileMode,FileName,FilePosition,fileReadable,fileTime,fileWritable,fillMatrix,findFiles,findHeft,FindOne,findProgram,findSynonyms,FiniteFittingIdeals,First,first,firstkey,FirstPackage,fittingIdeal,flagLookup,FlatMonoid,flatten,flattenRing,Flexible,flip,floor,flush,fold,FollowLinks,for,forceGB,fork,FormalGroupLaws,Format,format,formation,FourierMotzkin,FourTiTwo,fpLLL,frac,fraction,FractionField,frames,FrobeniusThresholds,from,fromDividedPowers,fromDual,Function,FunctionApplication,FunctionBody,functionBody,FunctionClosure,FunctionFieldDesingularization,fusePairs,futureParser,GaloisField,gb,GBDegrees,gbRemove,gbSnapshot,gbTrace,gcd,gcdCoefficients,gcdLLL,GCstats,genera,GeneralOrderedMonoid,GenerateAssertions,generateAssertions,generator,generators,Generic,GenericInitialIdeal,genericMatrix,genericSkewMatrix,genericSymmetricMatrix,gens,genus,get,getc,getChangeMatrix,getenv,getGlobalSymbol,getNetFile,getNonUnit,getPrimeWithRootOfUnity,getSymbol,getWWW,GF,gfanInterface,Givens,GKMVarieties,GLex,Global,global,globalAssign,globalAssignFunction,GlobalAssignHook,globalAssignment,globalAssignmentHooks,GlobalDictionary,GlobalHookStore,globalReleaseFunction,GlobalReleaseHook,Gorenstein,GradedLieAlgebras,GradedModule,gradedModule,GradedModuleMap,gradedModuleMap,gramm,GraphicalModels,GraphicalModelsMLE,Graphics,graphIdeal,graphRing,Graphs,Grassmannian,GRevLex,GroebnerBasis,groebnerBasis,GroebnerBasisOptions,GroebnerStrata,GroebnerWalk,groupID,GroupLex,GroupRevLex,GTZ,Hadamard,handleInterrupts,HardDegreeLimit,hash,HashTable,hashTable,HEAD,HEADER1,HEADER2,HEADER3,HEADER4,HEADER5,HEADER6,HeaderType,Heading,Headline,Heft,heft,Height,height,help,Hermite,hermite,Hermitian,HH,hh,HigherCIOperators,HighestWeights,Hilbert,hilbertFunction,hilbertPolynomial,hilbertSeries,HodgeIntegrals,hold,Holder,Hom,homeDirectory,HomePage,Homogeneous,Homogeneous2,homogenize,homology,homomorphism,HomotopyLieAlgebra,hooks,horizontalJoin,HorizontalSpace,HR,HREF,HTML,html,httpHeaders,Hybrid,HyperplaneArrangements,Hypertext,hypertext,HypertextContainer,HypertextParagraph,icFracP,icFractions,icMap,icPIdeal,id,Ideal,ideal,idealizer,identity,if,IgnoreExampleErrors,ii,image,imaginaryPart,IMG,ImmutableType,importFrom,in,incomparable,Increment,independentSets,indeterminate,IndeterminateNumber,Index,index,indexComponents,IndexedVariable,IndexedVariableTable,indices,inducedMap,inducesWellDefinedMap,InexactField,InexactFieldFamily,InexactNumber,InfiniteNumber,infinity,info,InfoDirSection,infoHelp,Inhomogeneous,input,Inputs,insert,installAssignmentMethod,installedPackages,installHilbertFunction,installMethod,installMinprimes,installPackage,InstallPrefix,instance,instances,IntegralClosure,integralClosure,integrate,IntermediateMarkUpType,interpreterDepth,intersect,intersectInP,Intersection,intersection,interval,InvariantRing,inverse,InverseMethod,inversePermutation,Inverses,inverseSystem,InverseSystems,Invertible,InvolutiveBases,irreducibleCharacteristicSeries,irreducibleDecomposition,isAffineRing,isANumber,isBorel,isCanceled,isCommutative,isConstant,isDirectory,isDirectSum,isEmpty,isField,isFinite,isFinitePrimeField,isFreeModule,isGlobalSymbol,isHomogeneous,isIdeal,isInfinite,isInjective,isInputFile,isIsomorphism,isLinearType,isListener,isLLL,isMember,isModule,isMonomialIdeal,isNormal,isOpen,isOutputFile,isPolynomialRing,isPrimary,isPrime,isPrimitive,isPseudoprime,isQuotientModule,isQuotientOf,isQuotientRing,isReady,isReal,isReduction,isRegularFile,isRing,isSkewCommutative,isSorted,isSquareFree,isStandardGradedPolynomialRing,isSubmodule,isSubquotient,isSubset,isSupportedInZeroLocus,isSurjective,isTable,isUnit,isWellDefined,isWeylAlgebra,ITALIC,Iterate,Jacobian,jacobian,jacobianDual,Jets,Join,join,Jupyter,K3Carpets,K3Surfaces,Keep,KeepFiles,KeepZeroes,ker,kernel,kernelLLL,kernelOfLocalization,Key,keys,Keyword,Keywords,kill,koszul,Kronecker,KustinMiller,LABEL,last,lastMatch,LATER,LatticePolytopes,Layout,lcm,leadCoefficient,leadComponent,leadMonomial,leadTerm,Left,left,length,LengthLimit,letterParser,Lex,LexIdeals,LI,Licenses,LieTypes,lift,liftable,Limit,limitFiles,limitProcesses,Linear,LinearAlgebra,LinearTruncations,lineNumber,lines,LINK,linkFile,List,list,listForm,listLocalSymbols,listSymbols,listUserSymbols,LITERAL,LLL,LLLBases,lngamma,load,loadDepth,LoadDocumentation,loadedFiles,loadedPackages,loadPackage,Local,local,localDictionaries,LocalDictionary,localize,LocalRings,locate,log,log1p,LongPolynomial,lookup,lookupCount,LowerBound,LUdecomposition,M0nbar,M2CODE,Macaulay2Doc,makeDirectory,MakeDocumentation,makeDocumentTag,MakeHTML,MakeInfo,MakeLinks,makePackageIndex,MakePDF,makeS2,Manipulator,map,MapExpression,MapleInterface,markedGB,Markov,MarkUpType,match,mathML,Matrix,matrix,MatrixExpression,Matroids,max,maxAllowableThreads,maxExponent,MaximalRank,maxPosition,MaxReductionCount,MCMApproximations,member,memoize,memoizeClear,memoizeValues,MENU,merge,mergePairs,META,method,MethodFunction,MethodFunctionBinary,MethodFunctionSingle,MethodFunctionWithOptions,methodOptions,methods,midpoint,min,minExponent,mingens,mingle,minimalBetti,MinimalGenerators,MinimalMatrix,minimalPresentation,minimalPresentationMap,minimalPresentationMapInv,MinimalPrimes,minimalPrimes,minimalReduction,Minimize,minimizeFilename,MinimumVersion,minors,minPosition,minPres,minprimes,Minus,minus,Miura,MixedMultiplicity,mkdir,mod,Module,module,ModuleDeformations,modulo,MonodromySolver,Monoid,monoid,MonoidElement,Monomial,MonomialAlgebras,monomialCurveIdeal,MonomialIdeal,monomialIdeal,MonomialIntegerPrograms,MonomialOrbits,MonomialOrder,Monomials,monomials,MonomialSize,monomialSubideal,moveFile,multidegree,multidoc,multigraded,MultigradedBettiTally,MultiGradedRationalMap,multiplicity,MultiplicitySequence,MultiplierIdeals,MultiplierIdealsDim2,MultiprojectiveVarieties,mutable,MutableHashTable,mutableIdentity,MutableList,MutableMatrix,mutableMatrix,NAGtypes,Name,nanosleep,Nauty,NautyGraphs,NCAlgebra,NCLex,needs,needsPackage,Net,net,NetFile,netList,new,newClass,newCoordinateSystem,NewFromMethod,newline,NewMethod,newNetFile,NewOfFromMethod,NewOfMethod,newPackage,newRing,nextkey,nextPrime,nil,NNParser,NoetherianOperators,NoetherNormalization,NonminimalComplexes,nonspaceAnalyzer,NoPrint,norm,normalCone,Normaliz,NormalToricVarieties,not,Nothing,notify,notImplemented,NTL,null,nullaryMethods,nullhomotopy,nullParser,nullSpace,Number,number,NumberedVerticalList,numcols,numColumns,numerator,numeric,NumericalAlgebraicGeometry,NumericalCertification,NumericalImplicitization,NumericalLinearAlgebra,NumericalSchubertCalculus,numericInterval,NumericSolutions,numgens,numRows,numrows,odd,oeis,of,ofClass,OL,OldPolyhedra,OldToricVectorBundles,on,OneExpression,OnlineLookup,OO,oo,ooo,oooo,openDatabase,openDatabaseOut,openFiles,openIn,openInOut,openListener,OpenMath,openOut,openOutAppend,operatorAttributes,Option,OptionalComponentsPresent,optionalSignParser,Options,options,OptionTable,optP,or,Order,order,OrderedMonoid,orP,OutputDictionary,Outputs,override,pack,Package,package,PackageCitations,PackageDictionary,PackageExports,PackageImports,PackageTemplate,packageTemplate,pad,pager,PairLimit,pairs,PairsRemaining,PARA,Parametrization,parent,Parenthesize,Parser,Parsing,part,Partition,partition,partitions,parts,path,pdim,peek,PencilsOfQuadrics,Permanents,permanents,permutations,pfaffians,PHCpack,PhylogeneticTrees,pi,PieriMaps,pivots,PlaneCurveSingularities,plus,poincare,poincareN,Points,polarize,poly,Polyhedra,Polymake,PolynomialRing,Posets,Position,position,positions,PositivityToricBundles,POSIX,Postfix,Power,power,powermod,PRE,Precision,precision,Prefix,prefixDirectory,prefixPath,preimage,prepend,presentation,pretty,primaryComponent,PrimaryDecomposition,primaryDecomposition,PrimaryTag,PrimitiveElement,Print,print,printerr,printingAccuracy,printingLeadLimit,printingPrecision,printingSeparator,printingTimeLimit,printingTrailLimit,printString,printWidth,processID,Product,product,ProductOrder,profile,profileSummary,Program,programPaths,ProgramRun,Proj,Projective,ProjectiveHilbertPolynomial,projectiveHilbertPolynomial,ProjectiveVariety,promote,protect,Prune,prune,PruneComplex,pruningMap,Pseudocode,pseudocode,pseudoRemainder,Pullback,PushForward,pushForward,Python,QQ,QQParser,QRDecomposition,QthPower,Quasidegrees,QuaternaryQuartics,QuillenSuslin,quit,Quotient,quotient,quotientRemainder,QuotientRing,Radical,radical,RadicalCodim1,radicalContainment,RaiseError,random,RandomCanonicalCurves,RandomComplexes,RandomCurves,RandomCurvesOverVerySmallFiniteFields,RandomGenus14Curves,RandomIdeals,randomKRationalPoint,RandomMonomialIdeals,randomMutableMatrix,RandomObjects,RandomPlaneCurves,RandomPoints,RandomSpaceCurves,Range,rank,RationalMaps,RationalPoints,RationalPoints2,ReactionNetworks,read,readDirectory,readlink,readPackage,RealField,RealFP,realPart,realpath,RealQP,RealQP1,RealRoots,RealRR,RealXD,recursionDepth,recursionLimit,Reduce,reducedRowEchelonForm,reduceHilbert,reductionNumber,ReesAlgebra,reesAlgebra,reesAlgebraIdeal,reesIdeal,References,ReflexivePolytopesDB,regex,regexQuote,registerFinalizer,regSeqInIdeal,Regularity,regularity,relations,RelativeCanonicalResolution,relativizeFilename,Reload,remainder,RemakeAllDocumentation,remove,removeDirectory,removeFile,removeLowestDimension,reorganize,replace,RerunExamples,res,reshape,ResidualIntersections,ResLengthThree,Resolution,resolution,ResolutionsOfStanleyReisnerRings,restart,Result,resultant,Resultants,return,returnCode,Reverse,reverse,RevLex,Right,right,Ring,ring,RingElement,RingFamily,ringFromFractions,RingMap,rootPath,roots,rootURI,rotate,round,rowAdd,RowExpression,rowMult,rowPermute,rowRankProfile,rowSwap,RR,RRi,rsort,run,RunDirectory,RunExamples,RunExternalM2,runHooks,runLengthEncode,runProgram,same,saturate,Saturation,scan,scanKeys,scanLines,scanPairs,scanValues,schedule,schreyerOrder,Schubert,Schubert2,SchurComplexes,SchurFunctors,SchurRings,SCRIPT,scriptCommandLine,ScriptedFunctor,SCSCP,searchPath,sec,sech,SectionRing,SeeAlso,seeParsing,SegreClasses,select,selectInSubring,selectVariables,SelfInitializingType,SemidefiniteProgramming,Seminormalization,separate,SeparateExec,separateRegexp,Sequence,sequence,Serialization,serialNumber,Set,set,setEcho,setGroupID,setIOExclusive,setIOSynchronized,setIOUnSynchronized,setRandomSeed,setup,setupEmacs,sheaf,SheafExpression,sheafExt,sheafHom,SheafOfRings,shield,ShimoyamaYokoyama,short,show,showClassStructure,showHtml,showStructure,showTex,showUserStructure,SimpleDoc,simpleDocFrob,SimplicialComplexes,SimplicialDecomposability,SimplicialPosets,SimplifyFractions,sin,singularLocus,sinh,size,size2,SizeLimit,SkewCommutative,SlackIdeals,sleep,SLnEquivariantMatrices,SLPexpressions,SMALL,smithNormalForm,solve,someTerms,Sort,sort,sortColumns,SortStrategy,source,SourceCode,SourceRing,SPACE,SpaceCurves,SPAN,span,SparseMonomialVectorExpression,SparseResultants,SparseVectorExpression,Spec,SpechtModule,SpecialFanoFourfolds,specialFiber,specialFiberIdeal,SpectralSequences,splice,splitWWW,sqrt,SRdeformations,stack,stacksProject,Standard,standardForm,standardPairs,StartWithOneMinor,stashValue,StatePolytope,StatGraphs,status,stderr,stdio,step,StopBeforeComputation,stopIfError,StopWithMinimalGenerators,Strategy,String,STRONG,StronglyStableIdeals,STYLE,Style,style,SUB,sub,SubalgebraBases,sublists,submatrix,submatrixByDegrees,Subnodes,subquotient,SubringLimit,Subscript,subscript,SUBSECTION,subsets,substitute,substring,subtable,Sugarless,Sum,sum,SumOfTwists,SumsOfSquares,SUP,super,SuperLinearAlgebra,Superscript,superscript,support,SVD,SVDComplexes,switch,SwitchingFields,sylvesterMatrix,Symbol,symbol,SymbolBody,symbolBody,SymbolicPowers,symlinkDirectory,symlinkFile,symmetricAlgebra,symmetricAlgebraIdeal,symmetricKernel,SymmetricPolynomials,symmetricPower,synonym,SYNOPSIS,syz,Syzygies,SyzygyLimit,SyzygyMatrix,SyzygyRows,syzygyScheme,TABLE,Table,table,take,Tally,tally,tan,TangentCone,tangentCone,tangentSheaf,tanh,target,Task,taskResult,TateOnProducts,TD,temporaryFileName,tensor,tensorAssociativity,TensorComplexes,terminalParser,terms,TEST,Test,testExample,testHunekeQuestion,TestIdeals,TestInput,tests,TEX,tex,TeXmacs,texMath,Text,TH,then,Thing,ThinSincereQuivers,ThreadedGB,threadVariable,Threshold,throw,Time,time,times,timing,TITLE,TO,to,TO2,toAbsolutePath,toCC,toDividedPowers,toDual,toExternalString,toField,TOH,toList,toLower,top,top,topCoefficients,Topcom,topComponents,topLevelMode,Tor,TorAlgebra,Toric,ToricInvariants,ToricTopology,ToricVectorBundles,toRR,toRRi,toSequence,toString,TotalPairs,toUpper,TR,trace,transpose,TriangularSets,Tries,Trim,trim,Triplets,Tropical,true,Truncate,truncate,truncateOutput,Truncations,try,TSpreadIdeals,TT,tutorial,Type,TypicalValue,typicalValues,UL,ultimate,unbag,uncurry,Undo,undocumented,uniform,uninstallAllPackages,uninstallPackage,Unique,unique,Units,Unmixed,unsequence,unstack,Up,UpdateOnly,UpperTriangular,URL,urlEncode,Usage,use,UseCachedExampleOutput,UseHilbertFunction,UserMode,userSymbols,UseSyzygies,utf8,utf8check,validate,value,values,Variable,VariableBaseName,Variables,Variety,variety,vars,Vasconcelos,Vector,vector,VectorExpression,VectorFields,VectorGraphics,Verbose,Verbosity,Verify,VersalDeformations,versalEmbedding,Version,version,VerticalList,VerticalSpace,viewHelp,VirtualResolutions,VirtualTally,VisibleList,Visualize,wait,WebApp,wedgeProduct,weightRange,Weights,WeylAlgebra,WeylGroups,when,whichGm,while,width,wikipedia,Wrap,wrap,WrapperType,XML,xor,youngest,zero,ZeroExpression,zeta,ZZ,ZZParser,
makeGWClass,getDiagonalClass,makeDiagonalForm,getSignature,isAnisotropic,isIsotropic,getAnisotropicPart,getSumDecompositionString,getGlobalA1Degree,getLocalA1Degree,isIsomorphicForm,addGW,getSumDecomposition}
}
\lstalias{Macaulay2output}{Macaulay2}

\usepackage{pmboxdraw}
\usepackage{newunicodechar}
\newunicodechar{┐}{\textSFiii}
\newunicodechar{┘}{\textSFiv}
\newunicodechar{└}{\textSFii}
\title{$C_p$-Mackey Functors in Macaulay2}
\author[T.\ Brazelton]{Thomas Brazelton}
\author[D.\ Chan]{David Chan}
\author[B.\ Mudrak]{Benjamin Mudrak}
\author[B.\ Spitz]{Ben Spitz}
\author[C.\ Vogeli]{Chase Vogeli}
\author[C.\ Wang]{Chenglu Wang}
\author[M.R.\ Zeng]{Michael R. Zeng}
\author[S.\ Zotine]{Sasha Zotine}
\date{July 2025}

\begin{document}

\begin{abstract}
We introduce the \texttt{CpMackeyFunctors} package for Macaulay2, which allows for computations with Mackey functors over a cyclic group of prime order.
\end{abstract}

\maketitle

\section{Introduction}

Mackey functors are algebraic objects which encode abstract induction and restriction operations parameterized by subgroups of a finite group $G$. These first arose in representation theory as an axiomatic framework for induction theorems for representation rings \cite{Green71,Dress}, but have since found applications to a wide range of contexts involving finite group actions. Among other examples, Mackey functors appear in nature in the following structures: 
\begin{itemize}
    \item group (co)homology and Tate cohomology of $G$-modules,
    \item the Bredon equivariant (co)homology of $G$-spaces \cite[Chapter XIII]{Alaska},
    \item algebraic $K$-theory of group rings \cite[Chapter 11]{Oliver} or of rings with $G$-action \cite{Brazelton},
    \item class groups, Mordell--Weil groups, and Shafarevich--Tate groups associated to $G$-Galois extensions of number fields (where they are often known as \textit{modulations}) \cite{Neukirch-cohomology,BleyBoltje}, and
    \item Grothendieck--Witt rings of quadratic forms \cite{CalleGinnett}.
\end{itemize}
We highlight that many of these examples fall within the realm of equivariant stable homotopy theory, where it has become standard to define invariants valued in Mackey functors. The standard analogy is that Mackey functors play the role that abelian groups do in nonequivariant algebraic topology. Just as calculations in algebraic topology frequently necessitate homological algebra, the development of homological algebra for Mackey functors is expected to have computational consequences in the equivariant setting.

In more detail, Mackey functors for a fixed group $G$ form a closed symmetric monoidal abelian category, so it is possible to make sense of familiar notions such as projective resolutions, $\Tor$, and $\Ext$ therein. However, these structures are difficult to access in concrete calculations, and all such computations are currently done by hand.

In the package \texttt{CpMackeyFunctors}, we remedy this by implementing Mackey functors for cyclic groups of prime order and the aforementioned constructions. To our knowledge, this is the first computer package for working with Mackey functors, and we hope it will prove useful throughout equivariant mathematics.

\subsection{Acknowledgments}
This work was carried out at the UW-Madison Macaulay2 Workshop in summer 2025 (DMS-2508868). We are grateful to the NSF and the university for their support. The first-named author is supported by NSF DMS-230324.  The second-named author was partially supported by NSF grant DMS-2135960. The fifth-named author was partially supported by NSF grant DMS-2052977.

\section{Mackey functors and homomorphisms}

For a finite group $G$, there are many equivalent definitions of a $G$-Mackey functor. Perhaps the most concise, due to Lindner \cite{Lindner}, is stating that a Mackey functor is a coproduct-preserving functor from the category of spans of finite $G$-sets to abelian groups:
\begin{align*}
    \Mack_G := \Fun^\oplus(\Span(\Fin_G),\Ab).
\end{align*}
Since every finite $G$-set decomposes into a direct sum of transitive ones, it suffices to specify the values of a Mackey functor on transitive $G$-sets, which are all torsors for $G/H$ for some subgroup $H\le G$, well-defined up to conjugacy. This leads to a \emph{finite} amount of data. In particular when the group is cyclic of prime order, the definition is very concise.

\begin{definition} Let $p$ be a prime number. A $C_p$\emph{-Mackey functor} $M$ is the data of two abelian groups $M(C_p/e)$ (called the \emph{underlying module}) and $M(C_p/C_p)$ (called the \emph{fixed module}), together with restriction, transfer, and conjugation homomorphisms of the form
    \begin{align*}
        \res \colon M(C_p/C_p) & \to M(C_p/e)   \\
        \tr \colon M(C_p/e)    & \to M(C_p/C_p) \\
        \conj \colon M(C_p/e)  & \to M(C_p/e),
    \end{align*}
    subject to the following axioms:
    \begin{enumerate}
        \item $\conj \circ \res = \res$ and $\tr\circ \conj = \tr$
        \item $\conj$ is an automorphism of order dividing $p$
        \item For each $x\in M(C_p/e)$ we have
        \begin{align*}
            \res(\tr(x)) = \sum_{i=0}^{p-1}\conj^i(x).
        \end{align*}
    \end{enumerate}
\end{definition}

In this package, we encode a $C_p$-Mackey functor as a new type, called \texttt{CpMackeyFunctor}. It is a hash table encoding the data of the prime, the underlying and fixed modules, and the three homomorphisms (restriction, transfer, and conjugation). To create a new Mackey functor, the user can call \texttt{makeCpMackeyFunctor(p,R,T,C)}. This takes in four pieces of data, namely the prime $p$ and the three homomorphisms $\res$, $\tr$, and $\conj$. For example, we can construct a very simple Mackey functor as follows:

\begin{computation}\label{comp:basic-constructor} \,
\begin{lstlisting}[language=Macaulay2]
i1 : needsPackage "CpMackeyFunctors";
i2 : M:=makeCpMackeyFunctor(2,id_(ZZ^1),matrix({{2}}),id_(ZZ^1))
o2 =        Res : | 1 |
       1  --------------->   1  -┐ Conj : | 1 |
     ZZ  <---------------  ZZ   <┘
            Tr : | 2 |
o2 : CpMackeyFunctor
\end{lstlisting} 
\end{computation}

\begin{notation} The data of a $C_p$-Mackey functor is often concisely encoded in what is called a \emph{Lewis diagram}, which displays the five pieces of data in the following shape:
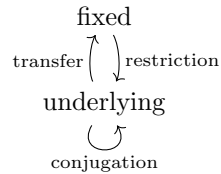
\[\begin{tikzcd}
	{\mathrm{fixed}} \\
	{\mathrm{underlying}}
	\arrow["{\text{restriction}}", curve={height=-6pt}, from=1-1, to=2-1]
	\arrow["{\text{transfer}}", curve={height=-6pt}, from=2-1, to=1-1]
	\arrow["{\mathrm{conjugation}}"', from=2-1, to=2-1, loop, in=300, out=240, distance=5mm]
\end{tikzcd}\]
\end{notation}

For example, the $C_2$-Mackey functor in \Cref{comp:basic-constructor} can be described by the following Lewis diagram
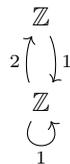
\[\begin{tikzcd}
	{\mathbb{Z}} \\
	{\mathbb{Z}}
	\arrow["1", curve={height=-6pt}, from=1-1, to=2-1]
	\arrow["2", curve={height=-6pt}, from=2-1, to=1-1]
	\arrow["1"', from=2-1, to=2-1, loop, in=300, out=240, distance=5mm]
\end{tikzcd}\]

\subsection{Examples of $C_p$-Mackey functors} 

We highlight some general examples of $C_p$-Mackey functors, as well as their method implementations in our package.

\begin{example}
    
The \emph{zero Mackey functor} is given by the following Lewis diagram:
\[\begin{tikzcd}
	0 \\
	0
	\arrow[curve={height=-6pt}, from=1-1, to=2-1]
	\arrow[curve={height=-6pt}, from=2-1, to=1-1]
	\arrow[from=2-1, to=2-1, loop, in=300, out=240, distance=5mm]
\end{tikzcd}\]
        This can be constructed as \texttt{makeZeroMackeyFunctor(p)} for any prime $p$.
\end{example}

\begin{example} There are a number of constructions of $C_p$-Mackey functors out of a single abelian group, possibly with an action by $C_p$.

\begin{enumerate} 
    
    \item Given any $C_p$-module $M$, we write $c\colon M\to M$ for the action of a chosen generator of $C_p$. We can write down the\emph{fixed point Mackey functor} of $M$, denoted $\FP(M)$, via the Lewis diagram:
    \[\begin{tikzcd}[ampersand replacement=\&]
    	{M^{C_p}} \\
    	M
    	\arrow["{\mathrm{incl.}}", curve={height=-6pt}, from=1-1, to=2-1]
    	\arrow["{\sum_{i=0}^{p-1} c^i}", curve={height=-6pt}, from=2-1, to=1-1]
    	\arrow["c"', from=2-1, to=2-1, loop, in=300, out=240, distance=5mm]
    \end{tikzcd}\]
    The restriction map is the inclusion, the transfer map is summation over orbits, and the conjugation map is $c\colon M\to M$. This is implemented as \texttt{makeFixedPointMackeyFunctor(p,c)}.

    \item Given any $C_p$-module $M$, we can define the \emph{orbits Mackey functor} via the Lewis diagram:
\[\begin{tikzcd}
	{M/C_p} \\
	M
	\arrow["{\sum_{i=0}^{p-1} c^i}", curve={height=-6pt}, from=1-1, to=2-1]
	\arrow["{\text{quot.}}", curve={height=-6pt}, from=2-1, to=1-1]
	\arrow["c"', from=2-1, to=2-1, loop, in=300, out=240, distance=5mm]
\end{tikzcd}\]
    This is implemented as \texttt{makeOrbitMackeyFunctor(p,c)}.

    \item Given a prime $p$ and an abelian group $M$, there is a $C_p$-Mackey functor that is zero on the underlying level, and $M$ on the fixed level:

\[\begin{tikzcd}
	M \\
	0
	\arrow["0", curve={height=-6pt}, from=1-1, to=2-1]
	\arrow["0", curve={height=-6pt}, from=2-1, to=1-1]
	\arrow["0"', from=2-1, to=2-1, loop, in=300, out=240, distance=5mm]
\end{tikzcd}\]
This can be constructed as \texttt{makeZeroOnUnderlyingMackeyFunctor(p,M)}.

\end{enumerate}

\end{example}

\begin{example}[Free Mackey functors] \label{eg:free-MFs}

There are \emph{free} Mackey functors which play a distinguished role in our construction of projective resolutions. We describe the sense in which they are ``free'' in \cref{prop:corepresentability}.

\begin{enumerate}

    \item \label{item:A-underline}
    The \textit{Burnside Mackey functor} for a finite group $G$ assigns to $G/H$ the Grothendieck group of finite $H$-sets. For $G=C_p$, this is given by the following Lewis diagram:
\[\begin{tikzcd}
	{\mathbb{Z}\oplus\mathbb{Z}t} \\
	{\mathbb{Z}}
	\arrow["{\text{res}}", curve={height=-6pt}, from=1-1, to=2-1]
	\arrow["t\cdot", curve={height=-6pt}, from=2-1, to=1-1]
	\arrow["1"', from=2-1, to=2-1, loop, in=300, out=240, distance=5mm]
\end{tikzcd}\]    
where $\mathbb{Z}\oplus\mathbb{Z}t \xto{\text{res}} \mathbb{Z}$ is the map $a + bt \mapsto a + bp$. For any prime $p$, this can be constructed as \texttt{makeBurnsideMackeyFunctor(p)} or \texttt{makeFixedFreeMackeyFunctor(p)}.

\item\label{item:B-underline}
We can also take the \emph{free $C_p$=Mackey functor on an underlying generator}, defined as
\[\begin{tikzcd}
	{\mathbb{Z}} \\
	{\mathbb{Z}[C_p]}
	\arrow[curve={height=-6pt}, from=1-1, to=2-1]
	\arrow[curve={height=-6pt}, from=2-1, to=1-1]
	\arrow["\gamma\cdot"', from=2-1, to=2-1, loop, in=300, out=240, distance=5mm]
\end{tikzcd}\]
where $\gamma$ is a generator for $C_p$. The restriction sends $1\mapsto \sum_{i=0}^{p-1} \gamma^i$, and the transfer sends $\gamma^i\mapsto 1$ for all $i$. This can be constructed via \texttt{makeUnderlyingFreeMackeyFunctor(p)}. We give this Mackey functor the special notation $\underline{B}$.

\end{enumerate}

\end{example}

\begin{example}[Examples from representation theory]

For a field $k$, there is a $G$-Mackey functor which assigns to $G/H$ the Grothendieck group of $H$-representations over $k$, where the transfer maps encode induction of representations and the restriction maps encode restriction of representations. Our package includes such examples for $G=C_p$ and $k=\mathbb R$ or $k=\mathbb C$.

\begin{enumerate}

\item The \textit{real representation $C_p$-Mackey functor} for $p$ odd has underlying module $\mathbb{Z}$ with trivial conjugation action, and fixed module $\mathbb{Z}\{\lambda_0,\lambda_1,\dots,\lambda_{(p-1)/2} \}$ where $\lambda_i$ is the two-dimensional real representation given by rotation by $(2\pi i)/p$ radians for $i > 0$ and $\lambda_0$ is the trivial one-dimensional representation. The restriction map is defined by $\lambda_i \mapsto 2$ for $i > 0$ and $\lambda_0 \mapsto 1$. The transfer map is defined by $x \mapsto x \cdot \sum_{0}^{(p-1)/2} \lambda_i$. Altogether, this $C_p$-Mackey functor is given by the following Lewis diagram:
\[\begin{tikzcd}
	{\mathbb{Z}\{\lambda_0,\lambda_1,\dots,\lambda_{(p-1)/2}   \}} \\
	{\mathbb{Z}}
	\arrow["\res", curve={height=-6pt}, from=1-1, to=2-1]
	\arrow["{x \mapsto x \cdot \sum_{0}^{(p-1)/2} \lambda_i}", curve={height=-6pt}, from=2-1, to=1-1]
	\arrow["1"', from=2-1, to=2-1, loop, in=300, out=240, distance=5mm]
\end{tikzcd}\]
If $p$ is even, then the real representation Mackey functor coincides with the complex representation Mackey functor as described in the next item.

This Mackey functor can be constructed as $\texttt{makeRealRepresentationMackeyFunctor(p)}$.

    \item The \textit{complex representation $C_p$-Mackey functor} has underlying module $\mathbb{Z}$ with trivial conjugation action, and fixed module $\mathbb{Z}\{\lambda_0,\lambda_1,\dots,\lambda_{p-1}\}$ where $\lambda_i$ is the one-dimensional complex representation given by multiplication by $e^{2\pi i/p}$. The restriction map is defined by $\lambda_i \mapsto 1$. The transfer map is defined by $x \mapsto x \cdot (\sum_0^{p-1}\lambda_i)$. Altogether, this $C_p$-Mackey functor is given by the following Lewis diagram:
\[\begin{tikzcd}
	{\mathbb{Z}\{\lambda_0,\lambda_1,\dots,\lambda_{p-1}   \}} \\
	{\mathbb{Z}}
	\arrow["{\lambda_i \mapsto 1}", curve={height=-6pt}, from=1-1, to=2-1]
	\arrow["{x \mapsto x \cdot \sum_{0}^{p-1} \lambda_i}", curve={height=-6pt}, from=2-1, to=1-1]
	\arrow["1"', from=2-1, to=2-1, loop, in=300, out=240, distance=5mm]
\end{tikzcd}\]
This Mackey functor can be constructed as $\texttt{makeComplexRepresentationMackeyFunctor(p)}$.

\end{enumerate}

\end{example}

\subsection{Mackey functor homomorphisms} 
\begin{definition}
A morphism of $C_p$-Mackey functors $f : M \to N$ consists of a pair of homomorphisms $(f_{C_p/e}, f_{C_p/C_p})$
\[\begin{tikzcd}[ampersand replacement=\&]
	{M(C_p/C_p)} \& {N(C_p/C_p)} \\
	{M(C_p/e)} \& {N(C_p/e)}
	\arrow["{f_{C_p/C_p}}", from=1-1, to=1-2]
	\arrow["\res", shift left=2, from=1-1, to=2-1]
	\arrow["\res", shift left=2, from=1-2, to=2-2]
	\arrow["\tr", shift left=2, from=2-1, to=1-1]
	\arrow["\conj"', from=2-1, to=2-1, loop, in=300, out=240, distance=5mm]
	\arrow["{f_{C_p/e}}"', from=2-1, to=2-2]
	\arrow["\tr", shift left=2, from=2-2, to=1-2]
	\arrow["\conj"', from=2-2, to=2-2, loop, in=300, out=240, distance=5mm]
\end{tikzcd}\]
such that for all $x\in M(C_p/e)$ and $y\in M(C_p/C_p)$ we have
\begin{enumerate}
    \item $f_{C_p/e}(\conj(x)) = \conj(f_{C_p/e}(x))$;
    \item $f_{C_p/e}(\res(y)) = \res(f_{C_p/C_p}(y))$;
    \item $f_{C_p/C_p}(\tr(x)) = \tr(f_{C_p/e}(x))$.
\end{enumerate}
In other words, $f$ commutes with $\res$, $\tr$, and $\conj$.
\end{definition}

We implement a type \texttt{MackeyFunctorHomomorphism}, which is again a hash table. It encodes the homomorphisms on the underlying and the fixed modules, and the source and target Mackey functors can be extracted.

\begin{example}
    If $M$ is a Mackey functor, we obtain its \emph{identity homomorphism} as \texttt{id\_M}.
\end{example}

\begin{example}
    If $f \colon M \to N$ and $g\colon N \to Q$ are composable Mackey functor homomorphisms, their composite can be obtained as \texttt{g*f}.
\end{example}


\subsection{The category of Mackey functors}

Mackey functors and Mackey functor homomorphisms form a category, denoted $\Mack_{C_p}$. By sending a Mackey functor to its fixed or underlying module, we obtain natural functors to abelian groups. A crucial fact is that these functors are corepresentable.
\begin{proposition}\label{prop:corepresentability}
Let $p$ be a prime number.
\begin{enumerate}
    \item The functor $\Mack_{C_p} \to \Ab$ sending $M$ to its fixed module $M(C_p/C_p)$ is corepresented by the \textit{Burnside Mackey functor} $\underline{A}$ (\Cref{eg:free-MFs}(\ref{item:A-underline})), in the sense that this forgetful functor is naturally isomorphic to the functor
    \[
        \Hom_{\Mack_{C_p}}(\underline{A},-) \colon \Mack_{C_p} \to \Ab.
    \]
    \item The functor $\Mack_{C_p}\to \Ab$ sending $M$ to its underlying module $M(C_p/e)$ is corepresented by the \textit{free Mackey functor on an underlying generator} $\underline{B}$ (\Cref{eg:free-MFs}(\ref{item:B-underline})), in the sense that this forgetful functor is naturally isomorphic to the functor
    \[
        \Hom_{\Mack_{C_p}}(\underline{B},-) \colon \Mack_{C_p} \to \Ab.
    \]
\end{enumerate}
\end{proposition}
Since epimorphisms in Mackey functors are exactly levelwise epimorphisms, this implies that both $\underline{A}$ and $\underline{B}$ are projective in the category $\Mack_{C_p}$. Moreover, the Yoneda lemma implies that every Mackey functor admits a surjection from some (possibly infinite) direct sum of copies of $\underline{A}\oplus \underline{B}$.  Thus $\underline{A}\oplus \underline{B}$ is a projective generator for $\Mack_{C_p}$ -- this fact will be very helpful to us, both for constructing random Mackey functors and in construction resolutions.  Before getting to this, we discuss the abelian category structure on $\Mack_{C_p}$.

\section{Homological algebra of \texorpdfstring{$C_p$}{Cₚ}-Mackey functors}

\subsection{The abelian category of Mackey functors} \label{subsec: abelian category}
Given any two $C_p$-Mackey functor homomorphisms $f,g \colon M \to N$, we can \emph{add them} by simply adding the module homomorphisms on the fixed and underlying levels. This gives the set $\Hom_{\Mack_{C_p}}(M,N)$ the structure of an abelian group; in particular we might say that $\Mack_{C_p}$ is \emph{pre-additive}.
Additionally, we can take \emph{direct sums} of two Mackey functors. If $M$ and $N$ are $C_p$-Mackey functors, their sum (implemented as \texttt{M++N}) is defined as
   \[
\left(
\begin{tikzcd}[ampersand replacement=\&,baseline=-1em]
	{M(C_p/C_p)} \\
	{M(C_p/e)}
	\arrow["\res", shift left=2, from=1-1, to=2-1]
	\arrow["\tr", shift left=2, from=2-1, to=1-1]
	\arrow["\conj"', from=2-1, to=2-1, loop, in=300, out=240, distance=5mm]
\end{tikzcd}
\right)
\oplus
\left(
\begin{tikzcd}[ampersand replacement=\&,baseline=-1em]
	{N(C_p/C_p)} \\
	{N(C_p/e)}
	\arrow["\res", shift left=2, from=1-1, to=2-1]
	\arrow["\tr", shift left=2, from=2-1, to=1-1]
	\arrow["\conj"', from=2-1, to=2-1, loop, in=300, out=240, distance=5mm]
\end{tikzcd}
\right)
=
\begin{tikzcd}[ampersand replacement=\&,row sep=huge]
	{M(C_p/C_p) \oplus N(C_p/C_p)} \\
	{M(C_p/e) \oplus N(C_p/e)}
	\arrow["{\begin{bmatrix} \res & 0 \\ 0 & \res \end{bmatrix}}", shift left=2, from=1-1, to=2-1]
	\arrow["{\begin{bmatrix} \tr & 0 \\ 0 & \tr \end{bmatrix}}", shift left=2, from=2-1, to=1-1]
	\arrow["{\begin{bmatrix} \conj & 0 \\ 0 & \conj \end{bmatrix}}"', from=2-1, to=2-1, loop, in=300, out=240, distance=5mm]
\end{tikzcd}.
\]

Finally, we have a well-defined \emph{kernel} and \emph{cokernel} of any Mackey functor homomorphism $f\colon M \to  N$, implemented as \texttt{ker f} and \texttt{coker f}, respectively, and defined as:
\[\begin{tikzcd}
	& {\ker(f_{C_p/C_p})} && {\operatorname{coker}(f_{C_p/C_p})} \\
	{\ker f =} & {\ker(f_{C_p/e})} & {\operatorname{coker}f=} & {\operatorname{coker}(f_{C_p/e})}
	\arrow["\res", shift left=2, from=1-2, to=2-2]
	\arrow["\res", shift left=2, from=1-4, to=2-4]
	\arrow["\tr", shift left=2, from=2-2, to=1-2]
	\arrow["\conj"', from=2-2, to=2-2, loop, in=300, out=240, distance=5mm]
	\arrow["\res", shift left=2, from=2-4, to=1-4]
	\arrow["\conj"', from=2-4, to=2-4, loop, in=300, out=240, distance=5mm]
\end{tikzcd}\]

Altogether, $\Mack_{C_p}$ has the structure of an \emph{abelian category}. Moreover, we will soon see that it is a \emph{closed symmetric monoidal} abelian category, which will allow us to define Ext and Tor groups therein. Before doing so, we discuss how the structures discussed so far can be used to construct random Mackey functors.

\subsection{Random Mackey functors and homomorphisms}

Every $C_p$-Mackey functor $M$ can be written as the cokernel of some map between projective Mackey functors. In particular, we have an exact sequence
\[
    \underline{A}^{\oplus k_1}\oplus \underline{B}^{k_2} \to \underline{A}^{\oplus \ell_1}\oplus \underline{B}^{\ell_2} \to M \to 0.
\]
Therefore to construct a ``random finitely generated $C_p$-Mackey functor'', it suffices to pick random integers $(k_1,k_2,\ell_1,\ell_2)$, construct a random map of the form above, then take its cokernel. By corepresentability (\Cref{prop:corepresentability}) this is equivalent to picking $k_1$ random elements in $(\underline{A}^{\oplus \ell_1}\oplus\underline{B}^{\oplus \ell_2})(C_p/C_p)$ and $k_2$ random elements in $(\underline{A}^{\oplus \ell_1}\oplus\underline{B}^{\oplus \ell_2})(C_p/e)$. Since Macaulay2 already includes functionality for choosing random elements of finitely generated modules, this is easy to implement.

This strategy is implemented in our method \texttt{makeRandomCpMackeyFunctor}, which takes as input a prime number $p$ and creates a random $C_p$-Mackey functor. 

\begin{computation}\label{comp:rand-mackey-functor} \,
\begin{lstlisting}[language=Macaulay2]
i2 : rand1 = makeRandomCpMackeyFunctor 3
o2 =                     Res : | 62 |
      cokernel | 65 |  ----------------> cokernel | 2015 |  -┐ Conj : | 521 |
                      <----------------                     <┘
                         Tr : | -1 |
o2 : CpMackeyFunctor
i3 : rand2 = makeRandomCpMackeyFunctor 3
o3 =                      Res : | -86 |
      cokernel | 135 |  -----------------> cokernel | 1935 |  -┐ Conj : | -179 |
                       <-----------------                     <┘
                          Tr : | 57 |
o3 : CpMackeyFunctor
\end{lstlisting} 
\end{computation}

With a similar method, we can use \texttt{makeRandomMackeyFunctorHomomorphism} to make a random homomorphism between \texttt{rand1} and \texttt{rand2}. 

\begin{computation}\label{comp:rand-mackey-functor} \,
\begin{lstlisting}[language=Macaulay2]
i4 : f = makeRandomMackeyFunctorHomomorphism(rand1, rand2)
o4 =                       fix : | -27 |
       cokernel | 135 | <------------------   cokernel | 65 |
             ^ |                                   ^ | 
             | |                                   | | 
             | v                                   | v 
      cokernel | 1935 | <------------------ cokernel | 2015 |
             ^ |           und : | -774 |          ^ | 
             └-┘                                   └-┘ 
o4 : MackeyFunctorHomomorphism
\end{lstlisting} 
\end{computation}

\subsection{The symmetric monoidal structure}
The category $\Mack_{C_p}$ has a closed symmetric monoidal structure, defined by the \emph{box product}
\begin{align*}
    -\boxtimes - \colon \Mack_{C_p} \times \Mack_{C_p} \to \Mack_{C_p}.
\end{align*}

\begin{definition}
    Given two $C_p$-Mackey functors $M$ and $N$, their \emph{box product}, denoted $M\boxtimes N$, is defined as the Mackey functor
    \[\begin{tikzcd}[ampersand replacement=\&]
	{\Big((M(C_p/C_p) \otimes N(C_p/C_p)) \oplus (M(C_p/e) \otimes N(C_p/e))\Big)/{\sim}} \\
	{M(C_p/e) \otimes N(C_p/e)}
	\arrow["{[a \otimes b, c \otimes d] \mapsto \res(a) \otimes \res(b) + \sum_{i=0}^{p-1} \conj^i (c) \otimes \conj^i (d)}", shift left=2, from=1-1, to=2-1]
	\arrow["{x \otimes y \mapsto [0, x \otimes y]}", shift left=2, from=2-1, to=1-1]
	\arrow["{\conj \otimes \conj}"', from=2-1, to=2-1, loop, in=300, out=240, distance=5mm]
\end{tikzcd}\]
where $\sim$ is the congruence generated by
\[\begin{aligned}
    (\tr(x) \otimes b, 0) &\sim (0, x \otimes \res(b)) \\
    (a \otimes \tr(y), 0) &\sim (0, \res(a) \otimes y) \\
    (0, x \otimes y)      &\sim (0, \conj(x) \otimes \conj(y)).
\end{aligned}\]
\end{definition}

The box product of two Mackey functors can be computed as \texttt{boxProduct(M,N)} or as \texttt{M**N}, as in the following example. Sometimes the Mackey functors and homomorphisms we obtain can be quite large. The \texttt{prune} method returns a smaller isomorphic presentation of the $C_p$-Mackey functor or a homomorphism by computing the minimal length representatives of the three structural maps. 

\begin{computation}\label{comp:rand-mackey-functor} \,
\begin{lstlisting}[language=Macaulay2]
i5 : prune(rand1 ** rand2)
o5 =                              Res : | 2 |
      cokernel | 5 0 0 0 0 0 |  ----------------> cokernel | 5 0 |  -┐ Conj : | 1 |
                               <----------------                    <┘
                                  Tr : | -1 |
o5 : CpMackeyFunctor
\end{lstlisting} 
\end{computation}

\begin{proposition}[Properties of the box product] \,
\begin{enumerate}
    \item We have natural isomorphisms $\underline{A}\boxtimes M \cong M \cong M \boxtimes\underline{A}$ for any $M\in \Mack_{C_p}$. That is, the Burnside Mackey functor $\underline{A}$ is the \emph{unit} for the symmetric monoidal structure on $\Mack_{C_p}$ defined by the box product.
    \item For any $M,N\in \Mack_{C_p}$ we have an isomorphism $M\boxtimes N \cong N \boxtimes M$, natural in both $M$ and $N$, such that the composite $M\boxtimes N\cong N\boxtimes M \cong M\boxtimes N$ is the identity.
\end{enumerate}
\end{proposition}
More concisely, we often say that the box product defines a \emph{symmetric monoidal structure} on $\Mack_{C_p}$. 

It turns out that $\Mack_{C_p}$ is \emph{closed} symmetric monoidal in the sense that the box product participates in a ``tensor-hom'' adjunction.

\begin{definition}
    Given two $C_p$-Mackey functors $M$ and $N$, we define their \emph{internal hom}, denoted $[M,N]\in \Mack_{C_p}$, as the following Mackey functor:
    \[\begin{tikzcd}[ampersand replacement=\&]
	{\Hom_{\Mack_{C_p}}(M,N)} \\
	{\Hom(M(C_p/e),N(C_p/e))}
	\arrow["{f \mapsto f_{C_p/e}}", shift left=2, from=1-1, to=2-1]
	\arrow["{h \mapsto (\sum_{i=0}^{p-1} \conj^i(h), \tr \circ h \circ \res)}", shift left=2, from=2-1, to=1-1]
	\arrow["\conj"', from=2-1, to=2-1, loop, in=300, out=240, distance=5mm]
\end{tikzcd}\]
\end{definition}
The statement that $\Mack_{C_p}$ is \emph{closed} symmetric monoidal is the statement that, for every $M,N,P\in \Mack_{C_p}$, we have an isomorphism
\[
    \Hom_{\Mack_{C_p}}(M \boxtimes N, P) \cong \Hom_{\Mack_{C_p}}(M, [N,P]),
\]
natural in $M$, $N$, and $P$.

The internal hom is implemented in the following method:

\begin{computation}\label{comp:rand-mackey-functor} \,
\begin{lstlisting}[language=Macaulay2]
i6 : prune internalHom(rand1,rand2)
o6 =                   Res : | -1 |
      cokernel | 5 |  ----------------> cokernel | 5 |  -┐ Conj : | 1 |
                     <----------------                  <┘
                        Tr : | 2 |
o6 : CpMackeyFunctor
\end{lstlisting} 
\end{computation}

\begin{proposition}[Properties of the internal hom] \,
\begin{enumerate}
    \item For any $M\in \Mack_{C_p}$, we have a natural isomorphism $[\underline{A},M]\cong M$.
    \item For any $M\in \Mack_{C_p}$, we have a natural isomorphism $[M,\underline{0}]\cong \underline 0$
    \item For any $X,Y\in \Mod_{\Z[C_p]}$, we have a natural isomorphism
    \[
        \FP \left( \Hom(X,Y) \right)\cong \left[\FP(X), \FP(Y)\right].
    \]
    where $\Hom(X,Y)$ denotes the \emph{internal} hom of $\mathbb{Z}[C_p]$-modules. 
\end{enumerate}
\end{proposition}

\subsection{Homological Algebra} Given the symmetric monoidal abelian category $\Mack_{C_p}$, we can now carry out homological algebra computations -- i.e. we can ask to take resolutions, derived functors, etc.

Our first method provides a \emph{resolution} of a $C_p$-Mackey functor $M$. This method takes in two inputs: $M$ and an integer $n$, and outputs a list of $C_p$-Mackey functor homomorphisms providing the first $n$ terms in a free resolution of $M$.

\begin{computation}\label{comp:resolution}\,
\footnotesize
\begin{lstlisting}[language=Macaulay2]
i2 : M:=makeCpMackeyFunctor(2,id_(ZZ^1),matrix({{2}}),id_(ZZ^1))
o2 =        Res : | 1 |
       1  --------------->   1  -┐ Conj : | 1 |
     ZZ  <---------------  ZZ   <┘
            Tr : | 2 |
o2 : CpMackeyFunctor
i3 : res(M,2)
o3 = {        fix : | 2 1 2 |       ,         fix : | -1 -1 -2 |       ,         fix : | -2 |          }
        1  <-------------------   3                 | 0  2  0  |                       | 0  |
      ZZ                        ZZ                  | 1  0  2  |                       | 1  |
      ^ |                       ^ |     3  <----------------------   3     3  <-------------------   1
      | |                       | |   ZZ                           ZZ    ZZ                        ZZ
      | v                       | v   ^ |                          ^ |   ^ |                       ^ | 
        1  <-------------------   3   | |                          | |   | |                       | | 
      ZZ      und : | 1 1 1 |   ZZ    | v                          | v   | v                       | v 
      ^ |                       ^ |     3  <----------------------   3     3  <-------------------   2
      └-┘                       └-┘   ZZ      und : | -1 0  -1 |   ZZ    ZZ      und : | -1 -1 |   ZZ
                                      ^ |           | 0  -1 -1 |   ^ |   ^ |           | -1 -1 |   ^ | 
                                      └-┘           | 1  1  2  |   └-┘   └-┘           | 1  1  |   └-┘ 
o3 : List
\end{lstlisting} 
\end{computation}

Since we have kernels and cokernels, we have access to (co)homology, and we can therefore compute derived functors. In particular we have implemented \texttt{Tor} and \texttt{Ext} computations. These methods take in two Mackey functors and an integer $i$, computing the $i$th Tor or Ext group.

\begin{computation}\label{comp:tor-ext}\,
\begin{lstlisting}[language=Macaulay2]
i4 : prune Tor_1(M,M)
o4 =                     Res : 0
     cokernel | 2 0 |  -----------> 0  -┐ Conj : 0
                      <-----------     <┘
                         Tr : 0
o4 : CpMackeyFunctor
i5 : prune Ext^4(M,M)
o5 =                   Res : 0
     cokernel | 2 |  -----------> 0  -┐ Conj : 0
                    <-----------     <┘
                       Tr : 0
o5 : CpMackeyFunctor
\end{lstlisting} 
\end{computation}

The following properties are immediate:

\begin{proposition}[Properties of Ext and Tor] For any $M,N\in \Mack_{C_p}$ we have
\begin{enumerate}
    \item $\Tor_0(M,N) = M \boxtimes N$
    \item $\Ext^0(M,N) = [M,N]$
    \item $\Ext^i(\underline{A},-)$ and $\Ext^i(\underline{B},-)$ are identically zero for $i>0$.
\end{enumerate}
\end{proposition}

We discuss an interesting conjecture arising from our computations of Ext-groups of $C_p$-Mackey functors in \Cref{sec:conjecture}, but we first discuss some examples of how to use this package to verify existing computations in the literature.

\section{Example Computations}

In this section we recall some computations with Mackey functors from the literature and demonstrate how they could have been performed using the \texttt{CpMackeyFunctors} package in Macaulay2.

\subsection{Equivariant algebraic \texorpdfstring{$K$}{K}-theory}

The first examples we consider are from a paper of Chan--Vogeli \cite{ChanVogeli} on the Galois-equivariant algebraic $K$-theory of finite fields. For any prime $p$, there is a unique degree $p$ extension $\mathbb{F}_{q^p}$ of the finite field $\mathbb F_q$. The algebraic $K$-groups of such fields were computed by Quillen.

\begin{theorem}[{\cite{Quillen}}]
    The algebraic $K$-groups of $\mathbb{F}_{q^p}$ are given by
    \[
        K_n(\mathbb{F}_{q^p}) \cong 
        \begin{cases}
            \mathbb{Z} & n=0\\
            \mathbb{Z}/(q^{ip}-1) & n = 2i-1\\
            0 & \mathrm{else}.
        \end{cases}
    \]
\end{theorem}
The Galois group of $\mathbb{F}_{q^p}$ over $\mathbb{F}_{q}$ is $G = C_p$, and the action of $C_p$ on $\mathbb{F}_{q^p}$ extends to an action on the algebraic $K$-groups. When $n = 2i-1>0$ is odd, we can pick a generator $\gamma$ of $C_p$ such that the action of $\gamma$ on $K_{2i-1}(\mathbb{F}_{q^p})\cong \mathbb{Z}/(q^{ip}-1)$ is multiplication by  $q^i$.  

Associated to the $C_p$-module $K_{2i-1}(\mathbb{F}_{q^p})$ is the orbit $C_p$-Mackey functor with Lewis diagram
\[\begin{tikzcd}[ampersand replacement=\&]
    	{\mathbb{Z}/(q-1)} \\
    	\mathbb{Z}/(q^{ip}-1)
    	\arrow["{\sum\limits_{j=0}^{p-1}q^{ij}}", curve={height=-6pt}, from=1-1, to=2-1]
    	\arrow["1", curve={height=-6pt}, from=2-1, to=1-1]
    	\arrow["q^i"', from=2-1, to=2-1, loop, in=300, out=240, distance=5mm]
    \end{tikzcd}\]
which, for brevity, we will denote by $R_i$. This Mackey functor is denoted by the symbol $\ominus^i$ in \cite{ChanVogeli}.

For any abelian group $M$ let us write $\langle M \rangle$ for the Mackey functor with $M$ as its fixed group and 0 as its underlying group.
The following result plays a crucial role in the resolution of extension problems arising from spectral sequences in \cite{ChanVogeli}.
\begin{proposition}[{\cite[Proposition 3.0.2]{ChanVogeli}}]
    For any primes $p$ and $q$, any abelian group $M$, and any $i$ we have $\mathrm{Ext}^1(R_i,\langle M\rangle) = \mathrm{Ext}^1(\langle M\rangle,R_i)=0$.
\end{proposition}

For any given abelian group and fixed values of $p$, $q$ and $i$ this computation can easily be implemented in the $C_p$-Mackey functors package for Macaulay2.  For instance, the following code snippet checks the proposition with the choices $M = \mathbb{Z}$, $p=5$, $q=7$, and $i=3$.

\begin{computation}\label{code: Chan Vogeli example} \,
\footnotesize
\begin{lstlisting}[language=Macaulay2]
i2 : p=5; q=7; i=2;
i5 : K = cokernel matrix {{q^(i*p)-1}}; -- the underlying level of Ri
i6 : gamma = map(K,K,q^i); -- action of Cp on K
o6 : Matrix K <-- K
i7 : Ri = makeOrbitMackeyFunctor (p,gamma) 
o7 =                               Res : | 5884901 |
     cokernel | -48 282475248 |  ---------------------> cokernel | 282475248 |  -┐ Conj : | 49 |
                                <---------------------                          <┘
                                   Tr : | 1 |
o7 : CpMackeyFunctor
i8 : M = makeZeroOnUnderlyingMackeyFunctor (p,ZZ^1) 
o8 =        Res : 0
       1  -----------> 0  -┐ Conj : 0
     ZZ  <-----------     <┘
            Tr : 0
o8 : CpMackeyFunctor
i9 : prune (Ext^1(Ri,M)).Fixed 
o9 = 0
o9 : ZZ-module
i10 : prune (Ext^1(M,Ri)).Fixed
o10 = 0
o10 : ZZ-module
\end{lstlisting} 
\end{computation}

\subsection{Equivariant cohomology}

Our next set of examples recovers calculations first preformed by Mingcong Zeng \cite{zeng_mackey_2018} in his work on equivariant cohomology. Before stating the results we need a definition.
\begin{definition}
    A $C_p$-Mackey functor is \emph{cohomological} if it has the property that $\tr\circ \res$ is multiplication by $p$. 
\end{definition}

The category $\mathrm{CohMack}_{C_p}$ of cohomological $C_p$-Mackey functors is an abelian subcategory of all $C_p$-Mackey functors, although it is not a Serre subcategory. If $M$ and $N$ are two cohomological $C_p$-Mackey functors we will write $\mathrm{ExtCoh}^i(M,N)$ for Ext computed in $\mathrm{CohMack}_{C_p}$.  Cohomological Ext is computed in the \texttt{CpMackeyFunctors} package using \texttt{ExtCoh(i,M,N)}.

\begin{theorem}[\cite{Arnold,BSW}]
    For any prime $p$ the global projective dimension of $\mathrm{CohMack}_{C_p}$ is $3$.
\end{theorem}

In particular, the groups $\mathrm{ExtCoh}^i(M,N)$ are zero for $i>3$. 

We now review the computational result of Zeng.  For a fixed prime $p$ we write $B_1$ for the $C_p$-Mackey functor $\langle \mathbb{Z}/p\rangle$.  We write $\underline{\mathbb{Z}}$ for the fixed point $C_p$-Mackey functor associated to the trivial $C_p$-action on $\mathbb{Z}$.  Both $B_1$ and $\underline{\mathbb{Z}}$ are cohomological.
\begin{proposition}[{\cite[Example 2.34]{zeng_mackey_2018}}]
    The cohomological Ext groups $\mathrm{ExtCoh}^i(B_1,\underline{\mathbb{Z}})$ are given by
    \[
        \mathrm{ExtCoh}^i(B_1,\underline{\mathbb{Z}}) =
        \begin{cases}
            B_1 & i=3\\
            0 & \mathrm{else}.
        \end{cases}
    \]
\end{proposition}

Since the global dimension of $\mathrm{CohMack}_{C_p}$ is $3$, the proposition is checked as soon as one checks that values of $\mathrm{ExtCoh}^i$ for $i\leq 3$.  Running the following code snippet confirms the proposition when $p=7$.

\begin{computation}\label{code: Zeng exmaple 1} \,
\footnotesize
\begin{lstlisting}[language=Macaulay2]
i2 : p=11;
i3 : B1 = makeZeroOnUnderlyingMackeyFunctor (p,cokernel matrix({{11}})) 
o3 =                    Res : 0
     cokernel | 11 |  -----------> 0  -┐ Conj : 0
                     <-----------     <┘
                        Tr : 0
o3 : CpMackeyFunctor
i4 : Z = makeFixedPointMackeyFunctor(p,id_(ZZ^1)) 
o4 =        Res : | 1 |
       1  ---------------->   1  -┐ Conj : | 1 |
     ZZ  <----------------  ZZ   <┘
            Tr : | 11 |
o4 : CpMackeyFunctor
i5 : for i to 3 do (
         print prune ExtCoh(i,B1,Z)
     )
     Res : 0
0  -----------> 0  -┐ Conj : 0
  <-----------     <┘
     Tr : 0
     Res : 0
0  -----------> 0  -┐ Conj : 0
  <-----------     <┘
     Tr : 0
     Res : 0
0  -----------> 0  -┐ Conj : 0
  <-----------     <┘
     Tr : 0
                   Res : 0
cokernel | 11 |  -----------> 0  -┐ Conj : 0
                <-----------     <┘
                   Tr : 0
\end{lstlisting} 
\end{computation}

For two cohomological $C_p$-Mackey functors $M$ and $N$ we can also consider the Tor cohomological Mackey functors $\mathrm{TorCoh}_i(M,N)$. Some interesting values of TorCoh groups are computed by Zeng.

\begin{proposition}[{\cite[Example 2.38]{zeng_mackey_2018}}]
    The cohomological Tor groups $\mathrm{TorCoh}_i(B_1,B_1)$ are given by
    \[
        \mathrm{TorCoh}_i(B_1,\underline{\mathbb{Z}}) =
        \begin{cases}
            B_1 & i=0,3\\
            0 & \mathrm{else}.
        \end{cases}
    \]
\end{proposition}

Cohomological Tor computations can easily be performed using the \texttt{CpMackeyFunctors} package.

\begin{computation}\label{code: Zeng exmaple 2} \,
\footnotesize
\begin{lstlisting}[language=Macaulay2]
i2 : p=11;
i3 : B1 = makeZeroOnUnderlyingMackeyFunctor (p,cokernel matrix({{11}})) 
o3 =                    Res : 0
     cokernel | 11 |  -----------> 0  -┐ Conj : 0
                     <-----------     <┘
                        Tr : 0
o3 : CpMackeyFunctor
i4 : for i to 3 do (
         print prune TorCoh(i,B1,B1)
     )
                       Res : 0
cokernel | 11 0 0 |  -----------> 0  -┐ Conj : 0
                    <-----------     <┘
                       Tr : 0
     Res : 0
0  -----------> 0  -┐ Conj : 0
  <-----------     <┘
     Tr : 0
     Res : 0
0  -----------> 0  -┐ Conj : 0
  <-----------     <┘
     Tr : 0
                                             Res : 0
cokernel | 11 0 0 0 0 0 0 0 0 0 0 0 0 0 |  -----------> 0  -┐ Conj : 0
                                          <-----------     <┘
                                             Tr : 0
\end{lstlisting} 
\end{computation}

\section{An interesting conjecture}\label{sec:conjecture}

It is well-known that the homological algebra of $C_p$-modules exhibits periodic behavior. In this section, we conjecture that $C_p$-Mackey functors exhibit similar behavior based on experiments run with our package $\texttt{CpMackeyFunctors}$. 

\subsection{Periodicity for $C_p$-modules}

We begin by reviewing the situation for $C_p$-modules.

It is a standard fact that the cohomology of modules over finite cyclic groups 2-periodic. In terms of $\Ext$-groups, this means that for any $C_p$-module $M$,
\[
    \Ext_{C_p}^i(\mathbb Z,M) \cong \Ext_{C_p}^{i+2}(\mathbb Z,M), \quad i>0.
\]  
This fact follows, for instance, from the explicit 2-periodic resolution of $\mathbb Z$ as a $C_p$-module \cite[XII.7]{CartanEilenberg}
\begin{equation} \label{eqn:Zres} \begin{tikzcd}
    \mathbb Z & \mathbb ZC_p \arrow[l, "\varepsilon"'] & 
    \mathbb ZC_p \arrow[l, "\gamma-1"'] & 
    \mathbb ZC_p \arrow[l, "N"'] &
    \cdots \arrow[l],
\end{tikzcd} \end{equation}
where $\varepsilon: \mathbb ZC_p\to\mathbb Z$ is the augmentation, and $N=1+\gamma+\cdots+\gamma^{p-1}\in \mathbb ZC_p$ is the additive norm element. In fact, this resolution implies a more fundamental periodicity result for $C_p$-modules.

\begin{proposition} \label{prop:Cp-modulePeriodicity}
    Any $C_p$-module $M$ admits a projective resolution 
    \[ \begin{tikzcd}
        M & P_0 \arrow[l] & P_1 \arrow[l, "\partial_1"'] & P_2 \arrow[l, "\partial_2"'] & \cdots \arrow[l, "\partial_3"']
    \end{tikzcd} \]
    such that $\partial_i=\partial_{i+2}$ for all $i>1$.
\end{proposition}

We thank Dave Benson for sharing the proof of this fact.

\begin{proof}
    Let $Z_0$ be a first syzygy of $M$, that is, the kernel of a surjection $P_0\to M$ where $P_0$ is projective. Then, 
    \[ \begin{tikzcd}[column sep = large]
    M & P_0 \arrow[l] &
    Z_0 \otimes \mathbb ZC_p \arrow[l] & 
    Z_0 \otimes \mathbb ZC_p \arrow[l,"\id\otimes(\gamma-1)"'] & 
    Z_0 \otimes \mathbb ZC_p \arrow[l,"\id\otimes N"'] &
    \cdots \arrow[l]
    \end{tikzcd} \]
    is a resolution of $M$ satisfying the desired property, where the first differential $Z_0\otimes \mathbb ZC_p\to P_0$ is the composite
    \[ \begin{tikzcd}[column sep = large]
        Z_0\otimes \mathbb ZC_p \arrow[->>,r, "\id\otimes\varepsilon"] & Z_0 \otimes \mathbb Z \cong Z_0 \arrow[r, hook] & P_0.
    \end{tikzcd} \qedhere \]
\end{proof}

\subsection{Periodicity conjectures for $C_p$-Mackey functors}

Motivated by sample computations on 1000 randomly generated $C_p$-Mackey functors over various primes, we conjecture that an analogue of \cref{prop:Cp-modulePeriodicity} holds for $C_p$-Mackey functors.

\begin{conjecture} \label{conjecture}
    Every $C_p$-Mackey functor admits an eventually $4$-periodic projective resolution.
\end{conjecture}

We record a consequence of \cref{conjecture} which may be easier to approach.

\begin{conjecture}
    For all $C_p$-Mackey functors $N,M$, there exists an integer $n_0$ such that
    \[\Ext^{n+4}(M,N) \cong \Ext^n(M,N)\]
    and
    \[\Tor_{n+4}(M,N) \cong \Tor_n(M,N)\]
    for all $n \geq n_0$.
\end{conjecture}

These conjectures will be investigated in future work.

\bibliography{bibliography}
\bibliographystyle{amsalpha}
\end{document}